\documentclass[12pt, reqno]{amsart}
\usepackage{graphicx}
\usepackage{hyperref}
\usepackage[caption = false]{subfig}
\usepackage{amsthm}
\usepackage{amsmath}
\usepackage[numbers]{natbib}
\usepackage{amssymb}
\usepackage{mathrsfs}
\usepackage{mathtools}
\usepackage{enumerate}
\usepackage{amssymb}
\usepackage{wrapfig}
\usepackage{float}
\usepackage{comment}
\usepackage{graphicx}
\allowdisplaybreaks

\usepackage[twoside,paperwidth=200mm, paperheight=297mm, top=35mm, bottom=35mm, left=30mm, right=26mm, marginparsep=3mm, marginparwidth=40mm]{geometry}

\numberwithin{equation}{section}

\theoremstyle{plain}

\newtheorem{theorem}{Theorem}[section]
\newtheorem{corollary}[theorem]{Corollary}

\newtheorem{lemma}{Lemma}[section]
\newtheorem{prop}[theorem]{Proposition}
\theoremstyle{definition}

\theoremstyle{remark}
\newtheorem{remark}{Remark}[section]

\makeatother

\begin{document}
	
	\title[Generalized De-preferential Random Graphs]{Generalised De-Preferential Random Graphs}
%	\thanks{I am thankful to etc.}	

	\author[Antar Bandhyopadhyay]{Antar Bandyopadhyay}
    \address[Antar Bandyopadhyay]{Theoretical Statistics and Mathematics Unit \\
         Indian Statistical Institute, Delhi Centre \\ 
         7 S. J. S. Sansanwal Marg \\
         New Delhi 110016 \\
         INDIA}
\email{antar@isid.ac.in}
\author[Kunal Joshi]{Kunal Joshi}
\address[Kunal Joshi]{
         Indian Statistical Institute, Delhi Centre \\ 
         7 S. J. S. Sansanwal Marg \\
         New Delhi 110016 \\
         INDIA}
\email{kunspuns07081999@gmail.com}

	\maketitle	

\begin{abstract}	
	We consider some further generalizations of the novel random graph models as introduced by Bandyopadhyay and Sen \cite{BaSe2025} and find asymptotic for the degree of a fixed vertex and along with the asymptotic degree distribution. We show that in the \emph{case of the inverse power law} the order of these statistics is much slower than the case of the simple inverse function, which was considered in \cite{BaSe2025}. However, the results for the linear case remain exactly the same even after introducing a "shift" parameter. \\



\end{abstract}

% keywords can be removed

\section{Introduction}

\subsection{Background}
From simulating the spread of diseases to traffic flow optimization, complex networks appears in many real life problems. Moreover, there is always some randomness associated with occurance of almost every real life event. Example, students forming friend circles in a class, people getting affected by virus in an epidemic, citing of a paper by researchers etc. Therefore, Random Graphs models are useful tools to explain some properties that are observed in various real life complex networks.

One such model is Preferential attachment random graphs in which new vertices more probably connects to those existing vertices with high degrees, one such model was proposed by Albert and Barabasi, since which a lot of work has been done in this field.

But in this project we work on the opposite models, where where new vertices are more probable to be attached to vertices with less degree. This model was introduced by Antar and Sen, in this project we extend that work to more general models.

\subsection{Model}

We use the same framework as in \cite{BaSe2025}. In \cite{BaSe2025}, linear and inverse de-preferential models were introduced. Both the models produces a graph sequence denoted by $\{G_n\}_{n=2}^{\infty}$, at each time $n$ a new vertex $v_{n+1}$ with $m(\geq1)$ half edges $e_{n+1,1}, \cdots, e_{n+1,m}$ is attached to the graph. The edges are attached sequentially and randomly to existing vertices so for every time $n$ the sequence produces a graph with a set of $n$ vertices denoted by $V(G_n) = \{v_1, \cdots, v_n\}$ and $mn$ edges. In our model, $multiple$ $edges$, $loops$ and $self$ $loops$ are not allowed. The edges of the new vertex are attached sequentially and degrees of the existing vertices are updated during intermediate steps. We denote the degree of a vertex $v_i$, $i$ $=$ $1$, $\cdots$, $n$ by $d_i(n+1,k)$, $k$ $=$ $0$, $\cdots$, $m-1$, after $k$ half edges of vertex $v_{n+1}$ is attached to the graph, degree referring to bot in and out degree of a vertex. We will use $\{ \mathcal F_{n,k}:n \geq 2, k = 0, \cdots, m-1\}$ to denote the natural filtration associated with the graph sequence. Note that $d_i(n+1,0) = d_i(n)$. For $m = 1$, the natural filtration will simply be denoted by $\{ \mathcal F_n: n \geq 2 \}$ and the half edge of the vertex $v_{n+1}$ will be denoted by $e_{n+1}$. Finally, let $N_k(n)$ denote the number of vertices of degree $k$ in $G_n$. $P_k(n) = \frac{N_k(n)}{n}$ denotes the empirical proportion of vertices of degree $k$ in $G_n$.

\subsubsection{Linear De-Preferential Model with a shift parameter $\theta \geq 1$ and $1\geq\alpha>0$}

For $m=1,$

\begin{equation}
P(e_{n+1} = \{v_{n+1}, v_i\}|\mathcal F_n) \propto \left( \theta - \frac{\alpha d_i(n)}{2n-1} \right)\\
\end{equation}

\begin{equation}
\implies P(e_{n+1} = \{v_{n+1}, v_i\}|\mathcal F_n) = \frac{1}{n\theta-\alpha}\left( \theta - \frac{\alpha d_i(n)}{2n-1} \right)
\end{equation}

For $m>1$, $j=1,\cdots,n$ and $k=0,1,\cdots,m-1$,

\begin{equation}
    P(e_{n+1} = \{v_{n+1}, v_i\}|\mathcal F_{n+1},k) = \frac{1}{n\theta-\alpha}\left( \theta - \frac{\alpha d_i(n+1,k)}{k+ (2n-1)m} \right)
\end{equation}

In \cite{BaSe2025}, the authors took $\theta = \alpha = 1$

\subsubsection{Inverse Power Law De-Preferential with exponent $\alpha > 0$ and "shift" parameter 
$\delta > -1$ for $m=1$}

\begin{equation}
P(e_{n+1} = \{v_{n+1}, v_i\}|\mathcal F_n) \propto \frac{1}{(\delta + d_i(n))^\alpha}
\end{equation}

\begin{equation}
    \implies P(e_{n+1} = \{v_{n+1}, v_i\}|\mathcal F_n) = \frac{C_n}{(\delta + d_i(n))^\alpha},
\end{equation}
where, 
\[
{C_n}^{-1} = D_n = \sum_{i=1}^{n}\frac{1}{(\delta + d_i(n))^\alpha}
\]
It is worth to note here that the proportionality constant $C_n^{-1}$ = $D_n$ may depend on the parameters $\alpha > 0$ and $\theta \geq 0$

For $m>1$, $j=1,\cdots,n$, $k=0,\cdots,m-1$
\begin{equation}
    P(e_{n+1} = \{v_{n+1}, v_i\}|\mathcal F_{n+1},k) = \frac{C_{n+1,k}}{(\delta + d_j(n+1,k))^\alpha}
\end{equation}
where,
\[
C_{n+1,k}^{-1} = D_{n+1,k} = \sum_{j=1}^{n}\frac{1}{((\delta + d_j(n+1,k))^\alpha)}
\]
In \cite{BaSe2025} $\delta$ and $\alpha$ were taken to be $0$ and $1$ respectively.

\section{Main Results}

\subsection{Linear De-Preferential with parameter $\theta$}

\begin{theorem}\label{thm:lde}
Let $(G_n)_{n=1}^{\infty}$ be a sequence of random graphs following a Linear De-Preferential with parameter $\alpha$ and $\theta$. Let $v_i$, $i=1,\cdots,n$ be a fixed vertex. Then as $n \to \infty$,
\begin{enumerate}
\item[1.]
Growth of degree of a fixed vertex $v_i$ is given by
\begin{equation}
\frac{d_i(n)}{m \log n} \xrightarrow{P} 1
\end{equation}

\item[2.]
We have the following CLT result
\begin{equation}
    \frac{d_i(n)-m\log n}{\sqrt{m\log  n}} \xrightarrow{d} N(0,1)
\end{equation}

\item[3.]
\begin{equation}
P(e_{n+1} = \{ v_{n+1}, v_i\}, d_{v_i}(n) = k \} \to \frac{1}{2^k}
\end{equation}

\item[4.]
Let,
\begin{equation}
    P_k(n) := \frac{1}{n}\sum_{i=1}^{n}\mathbf{1}_{(d_i(n)=k)}
\end{equation}
be a fraction of vertices with degree k. Then $\exists\ C_1>0,$ such that, as $n\to \infty$
\begin{equation}
    P\left( \max_{k}|P_k(n)-p_k|>C_1\left( \frac{\log n}{n} + \sqrt{\frac{\log n}{n}} \right) \right) = o(1)
\end{equation}
where $p_k=\frac{1}{2^k}$, $k\geq1.$ We have this result only for $m=1$ case

\begin{remark}
Note that results are same as in case when $\theta = \alpha = 1$
\end{remark}

\end{enumerate}
\end{theorem}

\subsection{Inverse De-Preferential with parameter $\alpha(\leq 1)$ and $\delta$}

\subsubsection{Taking $m=1$}

\begin{theorem}\label{thm:ida}
    Let $(G_n)_{n=1}^\infty$ be a sequence of random graphs following an Inverse De-Preferential with parameter $\alpha$. Let $v_i$, $i=1,\cdots,n$ be fixed vertices. Then as $n \to \infty$,

    \begin{enumerate}
        \item [1.]
        Growth of degree of a fixed vertex $v_i$ is given by

        \begin{equation}
            \frac{d_i}{(\log n)^\frac{1}{1+\alpha}} \xrightarrow{a.s.} \left( \frac{1+\alpha}{\lambda ^*} \right)^\frac{1}{1+\alpha}
        \end{equation}

        \item[2.]
        Let,
        \begin{equation}
          P_k(n) := \frac{1}{n}\sum_{i=1}^{n}\mathbf{1}(d_i(n)=k),
        \end{equation}
        be a fraction of vertices with degree k. Then $\forall$ k$\geq$1,
        \begin{equation}
            P_k(n) \xrightarrow{a.s.} \frac{(\delta+k)^\alpha\lambda^*}{(\delta+k)^\alpha\lambda^*+1}\prod_{i=1}^{k-1}\frac{1}{(\delta+k)^\alpha\lambda^*+1} 
        \end{equation}

        \item[3.]

        \begin{equation}
            P(e_{n+1} = \{ v_{n+1}, v_i\}, d_{v_i}(n) = k \}) \to \prod_{i=1}^{k-1}\frac{1}{(\delta+k)^\alpha\lambda^*+1} 
        \end{equation}

    \end{enumerate}
\end{theorem}

\begin{corollary}\label{cor:ida}
\textit{The limiting empirical degree distribution of a sequence of random graphs following the inverse de-preferential attachment model with $m = 1$ has mean $2$, mode $1$ and its tail decays at a rate faster than exponential.}
\end{corollary}

\subsubsection{Taking $m>1$}
\begin{theorem}\label{thm:ida_2}
Let $(G_n)_{n=1}^{\infty}$ be a sequence of random graphs following the inverse de-preferential attachment model with $m > 1$. Then $\exists$ constants $0 < c < C$, such that, as $n \to \infty$
\begin{equation}
    P \left( c \leq \frac{d_i(n)}{m \sqrt{\log n}} \leq C \right) \to 1.
    \tag{2.13}
\end{equation}
\end{theorem}

\section{Proofs for the Linear De-Preferential Models with parameter $\alpha$ and $\theta$}

\subsection*{Proof of Theorem~\ref{thm:lde}}
We have the following recursive relation for $E[d_i(n)]$
\begin{align*}
    E[d_i(n+1)] &= E[E[d_i(n+1)|\mathcal{F_{n+1,m-1}}]]\\
                &= E[E[d_i(n+1,m-1)+d_i(n+1)-d_i(n+1,m-1)|\mathcal{F_{n+1,m-1}}]] \\
                &= E[d_i(n+1,m-1)] + E[E[d_i(n+1)-d_i(n+1,m-1)|\mathcal(F_{n+1,m-1})]] \\ 
                &= E[d_i(n+1,m-1)] + \frac{1}{n\theta-\alpha}\left(\theta-\frac{\alpha d_i(n+1,m-1)}{m-1+(2n-1)m}\right) \\
                &=\left(1-\frac{\alpha}{n\theta-\alpha}\frac{1}{m-1+(2n-1)m}\right)E[d_i(n+1,m-1)]+\frac{\theta}{n\theta-\alpha} \\
                &=\alpha_nE[d_i(n,m)]+\frac{\theta\beta_n}{\theta n-\alpha}
\end{align*}

where
\[
\alpha_n = \prod_{j=0}^{m-1}\left(1-\frac{\alpha}{n\theta-\alpha}\frac{1}{j+(2n-1)m} \right)
\]
\[
\beta_n = 1 + \sum_{k=1}^{m-1}\prod_{j=m-k}^{m-1}\left(1-\frac{\alpha}{n\theta-\alpha}\frac{1}{j+(2n-1)m} \right)
\]

Let $a_n=E[d_i(n)]$, $n\geq i$
$a_i = m$ w.p 1

We define for $n \geq i+1$

\[
\gamma_n = \prod_{k=i}{n-1}\alpha_k
\]

We define $\gamma_i = 1$

We have the recursion

\[
\frac{a_{m+1}}{\gamma_{n+1}} = \frac{a_n}{\gamma_n} + \frac{\theta\beta_n}{(\theta n - \alpha)\gamma_{n+1}}
\]

\[
\implies c_{n+1} = c_n + \frac{\theta\beta_n}{(\theta n - \alpha)\gamma_{n+1}}
\]
where $c_n = \frac{a_n}{\gamma_n}$, $n\geq i$

Therefore, we have 
\[
c_{n+1} = m + \frac{\theta\beta_i}{(\theta i -\alpha)\gamma_{i+1}} + \cdots + \frac{\theta\beta_n}{(\theta n -1)\gamma_{n+1}}
\]

$\beta_n\uparrow m$ as $n\to \infty$
and 
\[
\gamma_n \downarrow \kappa = \prod_{j=0}^{m-1}\prod_{n=i}^{\infty}\left(1-\frac{\alpha}{n\theta-\alpha}\frac{1}{j+(2n-1)m} \right)
\]
as $n\to \infty$

Let $h_n=1+\frac{1}{i+1}+\cdots+\frac{1}{n}$
then $\frac{h_n}{\log n}\to 1$ as $n \to \infty$
Using these results, we can prove that $\frac{c_n}{h_n}\to m\kappa^{-1}$ as $n\to \infty$, which further implies $\frac{a_n}{m\log n}\to1$ as $n\to \infty$

The above calculations were done to show that calculations done for case $\theta = \alpha = 1$ in \cite{BaSe2025} also holds for general case, which can be easily verified.

%%Proofs for inverse de-preferential case

\section{Proofs for Inverse De-preferential case with parameter $\alpha(\leq1)$ and $\delta$ }

In this section we provide proofs of the main results for the inverse de-preferential case.

\subsection{Embedding}
We use $Crump-Mode-Jagers (CMJ)$ branching process for $m=1$ case and $Athreya-Carlin$ $Embedding$ for $m>1$ case as done in \cite{BaSe2025} 

\subsubsection{CMJ Branching Process}\label{section:4.1.1}
$\mathcal{G}$ denotes the space of finite rooted trees. If $\mathcal{T} \in \mathcal{G}$ and $x$ is a vertex of $\mathcal{T}$, then we define $(\mathcal{T})\downarrow_{x}$ as the sub-tree consisting of the descendants of $x$.

We start with a graph consisting of a vertex with a half edge. Each vertex reproduces independently according to i.i.d.\ copies of a pure birth process $\{\xi(t): t \geq 0\}$. $\xi(0) = 1$ w.p. $1$ and
\[
P(\xi(t+h) = k+1 \mid \xi(t) = k) = \frac{h}{(k+\delta)^\alpha} + o(h)
\]

Let $\{\Upsilon(t) : t \geq 0\}$ denote the randomly growing tree. For every $t \geq 0$, $\Upsilon(t) \in \mathcal{G}$. We define the following sequence of random times
\begin{align*}
    \tau_2 &:= \inf\{t \geq 0 : |\Upsilon(t)| = 2\} \\
    \tau_3 &:= \inf\{t \geq \tau_2 : |\Upsilon(t)| = 3\} \\
    &\vdots \\
    \tau_n &:= \inf\{t \geq \tau_{n-1} : |\Upsilon(t)| = n\} \\
    &\vdots
\end{align*}
We look at the tree $\Upsilon(t)$ at the random times $\{\tau_n\}$. Let $\{G_n\}_{n=2}^\infty$ denote the random graph sequence under the inverse de-preferential attachment case when $m=1$. This gives us the following result.

\begin{theorem}\label{thm:4.1}
\textit{The sequence of random graphs $\{G_n\}_{n=2}^\infty$ is distributed identically as $\{\Upsilon(\tau_n)\}_{n=2}^\infty$.}
\end{theorem}

We will find it useful to consider the pure birth process $\{\xi(t): t \geq 0\}$ as a point process, where the points occur at the birth times $\{T_n\}_{n=1}^\infty$ of the pure birth process. We define the expected Laplace Transform of the point process $\{\xi(t): t \geq 0\}$
\[
\hat{\rho}(\lambda;\alpha,\delta) = E \left( \int_0^\infty \exp(-\lambda t) d\xi(t) \right)
    = \sum_{n=1}^\infty \prod_{i=0}^{n-1} \frac{1}{(i+\delta+1)^\alpha\lambda + 1}
\]
$\hat{\rho}(\lambda;\alpha,\delta)$ may be calculated easily because $\{T_n - T_{n-1}\}$ are independent random variables with $Exp(1/(\delta+n)^\alpha)$ distribution, that is, exponential distributions with mean $n$. We observe that the equation $\hat{\rho}(\lambda;\alpha,\delta) = 1$ has a unique root $\lambda = \lambda^*$. $\lambda^*$ is usually referred to as the Malthusian Parameter in the context of Crump-Mode-Jagers Branching Processes.

We don't have a closed formulae for $\lambda^*$ but we can study how it is dependent on $\delta$ and $\alpha$.

\begin{figure}[h]
  \centering
  \includegraphics[width=0.8\textwidth]{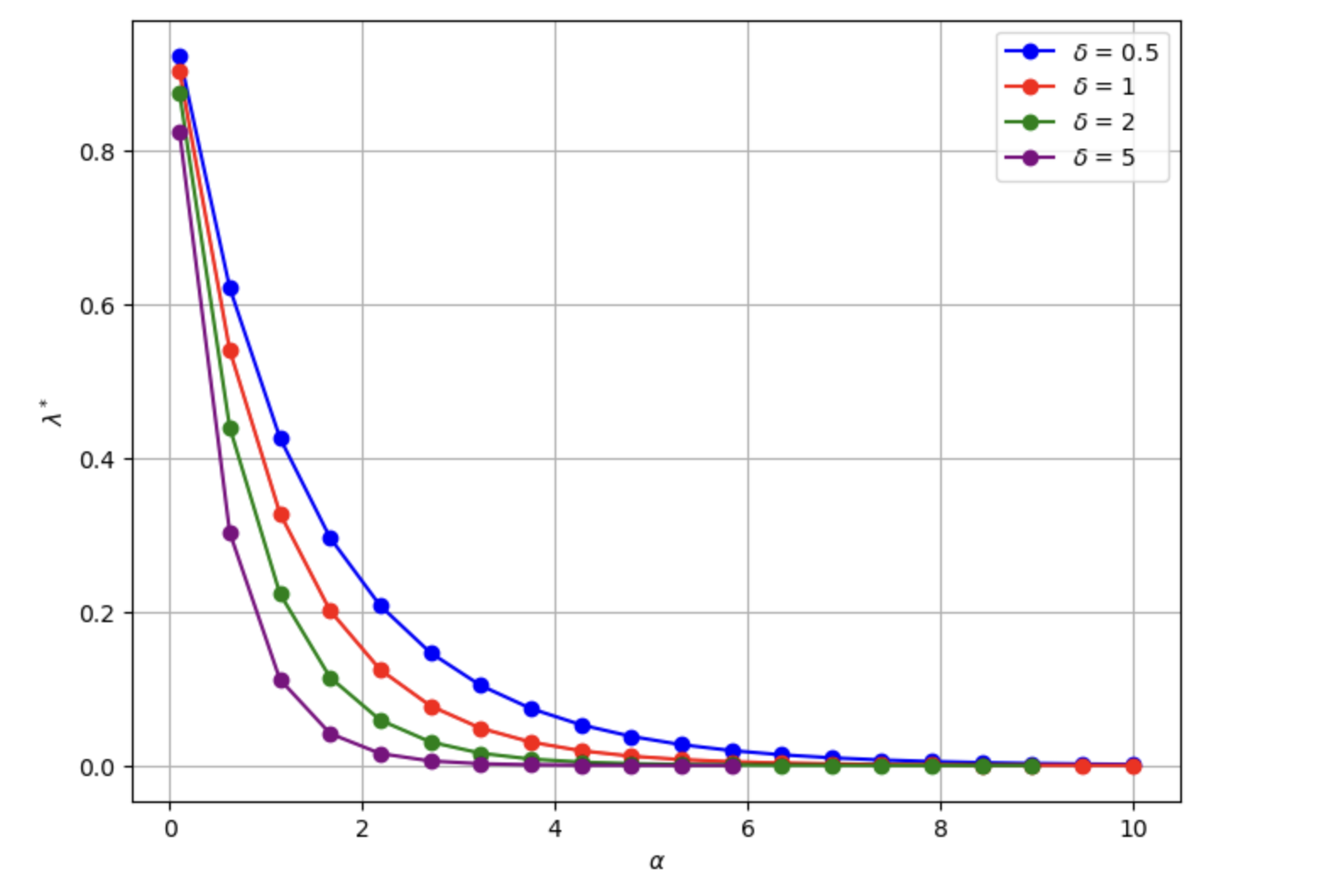}
  \caption{$\lambda^*$ vs $\alpha$ for multiple $\delta$ values}
  \label{fig:lambda-alpha}
\end{figure}

\begin{figure}[h]
  \centering
  \includegraphics[width=0.8\textwidth]{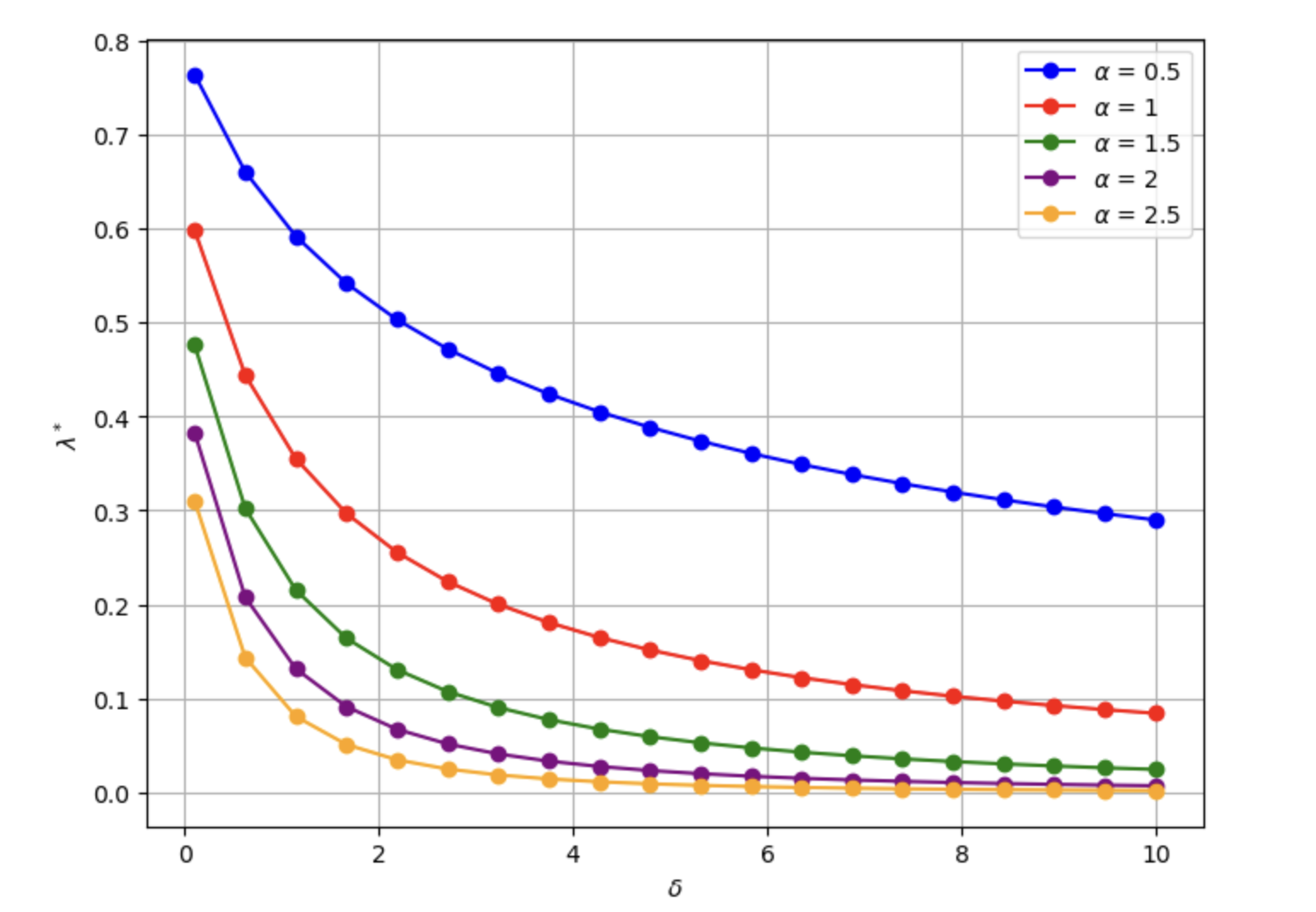}
  \caption{$\lambda^*$ vs $\delta$ for multiple $\alpha$ values}
  \label{fig:lambda-delta}
\end{figure}

From Figures 1 and 2,  appears to be a monotonically decreasing function of both $\alpha$ and $\delta$. Indeed, for fixed $\lambda>0$ and $\delta$, each factor $\frac{1}{(i+\delta+1)^\alpha\lambda + 1}$ is strictly ecreasing in $\alpha$, so the product and hence $\hat{\rho}(\lambda;\alpha,\delta)$ decrease as $\alpha$ increases. To keep the identity $\hat{\rho}(\lambda^*;\alpha,\delta)=1$ valid, $\lambda^*$ must therefore needs to be decreased. Similarly in case $\delta$ increases.

%This is because if we increase $\alpha$ and keep $\lambda^*$ fixed or decrease each factor $\frac{1}{(i+\delta+1)^\alpha\lambda + 1}$ of $\hat{\rho}(\lambda)$ will also decrease and hence won't add upto 1, therefore $\lambda^*$ needs to be increased and similarly it needs to be decreased if $\delta$ increases.

The process $\{\Upsilon(t) : t \geq 0\}$ is a supercritical, Malthusian Branching Process. Then using a theorem from O. Nerman (1981)(Theorem~A,~\cite{RudTotVal2007}), we have the following result.

\begin{theorem}\label{thm:4.2}
\textit{Consider a bounded function $\phi : \mathcal{G} \to \mathbb{R}$. Then the following limit holds almost surely}
\[
    \lim_{t \to \infty} \frac{1}{|\Upsilon(t)|} \sum_{x \in \Upsilon(t)} \phi(\Upsilon(t)\downarrow_{x}) = \lambda^* \int_0^\infty \exp\{-\lambda^* t\} E(\phi(\Upsilon(t))) dt.
\]
\end{theorem}

\subsubsection{Athreya-Karlin Embedding}\label{section:4.1.2}
For $m > 1$, we will couple our graph process with a sequence of Yule processes (i.i.d pure birth processes), with appropriate birth rates, such that, our degree sequences at each vertex will have the same distribution as the number of particles in the respective Yule processes sampled at suitable random times. Similar coupling has been used in~\cite{Ath2007},\cite{AthGhoSeth2008} and ~\cite{BuDa90}, where it is termed as \textit{Rubin's construction} in the context of reinforced random walks. The coupling with our specific birthrates, namely, $\lambda_i = \frac{1}{i}$ has also appeared in a recent work of Thacker and Volkov~\cite{ThVol2018}.

It is worth mentioning here that unlike in the previous case $m = 1$, when the entire graph process can be embedded inside a branching tree, in this case we are not embedding the entire graph process into the random forest which can be obtained from the Yule processes. Instead, this coupling is essentially only for the degree sequences at each vertex. Thus, having \textit{self-loops and multiple edges} in our random graphs, which can occur since $m > 1$, do not create any contradiction to this coupling.

Let $\{Z(t): t \geq 0\}$ be a pure birth process with $Z(0) = m$ w.p.\ 1 and birth rates $\{\lambda_i\}_{i=m}^\infty$, $\lambda_i = \frac{1}{(i+\delta)^\alpha}$. Let $\{Z_i(t): t \geq 0\}_{i=1}^\infty$ be i.i.d.\ copies of the pure birth process $Z(t)$.

We will define a sequence of random times $\{\tau_n\}_{n=1}^\infty$. Let $\tau_1 = 0$ w.p.1. We start the process $Z_1(t)$ at $t = 0$. Let $T_1^{(2)}$ be the time after $\tau_1$ when the first birth takes place in $Z_1$. Let $T_2^{(2)}$ be the time after $\tau_1 + T_1^{(2)}$ when the second birth occurs in $Z_1$. We continue in this manner to get $T_1^{(2)}, \cdots, T_m^{(2)}$. Let $\tau_2$ be the time when the $m^{\text{th}}$ birth occurs in $Z_1(t - \tau_1)$. We start a new process $Z_2(t)$ at $t = \tau_2$. We have,
\[
\tau_2 - \tau_1 = T_1^{(2)} + \cdots + T_m^{(2)}.
\]
In general, let $T_k^{(n+1)}$ denote the time of $k^{\text{th}}$ cumulative birth in $Z_1, \ldots, Z_n$ after $\tau_n + T_1^{(n+1)} + \cdots + T_{k-1}^{(n+1)}$. Let $\tau_{n+1}$ be the time after $\tau_n$ when the $m^{\text{th}}$ birth takes place in the processes $Z_1, Z_2, \ldots, Z_n$ after $\tau_n$. We start the process $Z_{n+1}(t)$ at $t = \tau_{n+1}$. Therefore, we have,
\[
\tau_{n+1} - \tau_n = T_1^{(n+1)} + \cdots + T_m^{(n+1)}
\]
We define, for $j=1, \cdots, n$,
\[
\tilde{d}_j(n+1, k) = Z_j(\tau_n + T_1^{(n+1)} + \cdots + T_k^{(n+1)} - \tau_j)
\]
We have the following embedding result.
\begin{theorem}\label{thm:4.3}
The sequence $\{\tilde{d}_j(n+1, k), k=0,\cdots, m-1, j=1, \cdots, n, n \geq 2\}$ and $\{d_j(n+1, k), k=0, \cdots, m-1, j=1, \cdots, n, n \geq 2\}$ are identically distributed.
\end{theorem}

It is important to emphasize at this point that the embedding of the random graph process in the $m = 1$ case into a CMJ branching process induces an Athreya-Karlin Embedding of the random graph process in the same probability space. We will be utilizing both these embeddings for the $m = 1$ case to establish certain properties of these random graphs.

We note that the pure birth process $Z(t)$ considered in Section~\ref{section:4.1.2} reduces to the birth process $\xi(t)$ considered in Section~\ref{section:4.1.1} if we fix $m=1$. We will need the following results about the asymptotic properties of these pure birth processes.

\begin{theorem}\label{thm:4.4}
Let $\{Z(t): t \geq 0\}$ be a pure birth process with $Z(0) = m$ w.p.\ 1 and birth rates $\{\lambda_i\}_{i=m}^\infty,\, \lambda_i = \frac{1}{(i+\delta)^\alpha}$. Then $\frac{Z(t)}{{t}^{\frac{1}{1+\alpha}}} \xrightarrow{P} (1+\alpha)^{\frac{1}{1+\alpha}}$ as $t \to \infty$.
\end{theorem}

\emph{Proof.} Let
\[
T_1 = 0 \qquad
T_2 = \inf\{t \geq T_1 : Z(t) = m+1\} \qquad
T_3 = \inf\{t \geq T_2 : Z(t) = m+2\}
\]
\[
\vdots
\]
Then we know, $T_1,\, T_2 - T_1,\, \cdots,\, T_n - T_{n-1},\, \cdots$ are independent exponential random variables. Therefore, $L_n = T_{n+1} - T_n$ are independent and $L_n \sim Exp(\frac{1}{(m+n-1+\delta)^\alpha})$. Then we have
\[
T_n = L_1 + \cdots + L_{n-1}
\implies E(T_n) = (m+\delta)^\alpha + (m+1+\delta)^\alpha + \cdots + (m+n-2+\delta)^\alpha
\]
Also, $\operatorname{Var}(L_n) = (m+n-1+\delta)^{2\alpha}$. Hence we have,
\[
\sum_{k=1}^{\infty} \operatorname{Var} \left(\frac{L_k}{k^{1+\alpha}}\right)
= \sum_{k=1}^{\infty} \frac{(m+k-1+\delta)^{1+\alpha}}{k^{2+2\alpha}} < \infty.
\]
Now observe $\sum\operatorname{Var}\left( \frac{L_k}{k^{1+\alpha}} \right) \sim \sum_k^\infty\frac{1}{k^{2}}<\infty$.
Thus, $\sum_{k=1}^{\infty} \frac{L_k - E[L_k]}{k^{1+\alpha}} < \infty$ w.p. $1$. First, this implies that $\frac{L_k}{k^2} \to 0$ a.s., because $E[L_k] = (m+k-1)$. Further, an application of Kronecker's Lemma yields that $\frac{T_n - E[T_n]}{n^{1+\alpha}} \to 0$ w.p. $1$. This allows us to conclude that $\frac{T_n}{n^{1+\alpha}} \to \frac{1}{1+\alpha}$ a.s.

We observe that $Z(t) \uparrow \infty$ as $t \to \infty$ and $T_{Z(t)} \leq t < T_{Z(t)+1}$. Therefore,
\[
\frac{T_{Z(t)}}{(Z(t))^{1+\alpha}} \leq \frac{t}{(Z(t))^{1+\alpha}} < \frac{T_{Z(t)}}{(Z(t))^{1+\alpha}} + \frac{L_{Z(t)}}{(Z(t))^{1+\alpha}}
\implies \frac{Z(t)}{t^{1+\alpha}} \to (1+\alpha)^{1+\alpha} \qquad a.s.
\]

This completes the proof. \qed

\subsection{Technical Results on Normalizing Constant}
The results established below allow us to approximate the normalizing constants for the de-preferential attachment model.

Let $\widetilde{\mathcal{F}}_{n,k}$ denote the natural filtration associated with the continuous time embedding $\tilde{d}_j(n+1, k)$ described in section~\ref{section:4.1.2}. We define
\[
\widetilde{C}_{n+1,k}^{-1} = \widetilde{D}_{n+1,k} = \sum_{j=1}^n \frac{1}{\tilde{d}_j(n+1, k)}.
\]
The natural filtration $\{\widetilde{\mathcal{F}}_n\}$ and the constants $\widetilde{C}_n$ and $\widetilde{D}_n$ are defined analogously for the $m=1$ case.

\begin{lemma}\label{lemma:4.1}
For $m=1$, $\frac{\tilde{D}_n}{n} \to \lambda^*$ a.s.
\end{lemma}

\emph{Proof.} The proof follows from Theorem~A, \cite{RudTotVal2007} by using the bounded functional $\phi(G) = \frac{1}{1+\# \text{children}(\emptyset)}$, where $\emptyset$ is the root of a finite tree $G \in \mathcal{G}$. \qed

For $m > 1$, we have the following bounds on the normalizing constants.

\begin{lemma}\label{lemma:4.2}
For all $n \geq 2$, for all $k = 0, 1, \ldots, m-1$, $\frac{(m+\delta)^\alpha}{n-1} \leq C_{n,k} \leq \frac{(2m+\delta)^\alpha}{n-1}$ w.p.\ 1.
\end{lemma}

\emph{Proof.} We note that $C^{-1}_{n,k} = D_{n,k} = \sum_{j=1}^{n-1} \frac{1}{(d_j(n, k) + \delta)^\alpha}$. We observe that $d_j(n, k) \geq m$ for $j=1, \ldots, n-1$ and therefore $D_{n,k} \leq \frac{n-1}{(m+\delta)^\alpha}$. We also observe that
\[
\sum_{j=1}^{n-1} d_j(n, k) = k + m(2n-3).
\]
Therefore, by the A.M.-H.M. inequality,
\[
\sum_{j=1}^{n-1} \frac{1}{d_j(n, k)}
\geq \frac{(n-1)^2}{\sum_{j=1}^{n-1}(d_j(n,k)+\delta)^\alpha}
\]

Using Jensen's inequality for summation we get
\begin{align}
\left( \sum_{j=1}^{n-1}(d_j(n,k)+\delta)^\alpha \right) &\leq n^{1-\alpha} \left( \sum_{j=1}^{n-1}(d_j(n,k)+\delta) \right)^\alpha\\
&=n^{1-\alpha}(k+m(2n-3)+(n-1)\delta)^\alpha\\
\implies\frac{(n-1)^2}{\sum_{j=1}^{n-1}(d_j(n,k)+\delta)^\alpha} &\geq\frac{(n-1)^{1+\alpha}}{(k+m(2n-3)+(n-1)\delta)^\alpha}\\
&\geq \frac{(n-1)^{1+\alpha}}{(k+m(2n-3)+(n-1)\delta + m-k)^\alpha}\\
&=\frac{n-1}{(2m+\delta)^\alpha}
\end{align}
Combining, we get that $\frac{n-1}{(2m+\delta)^\alpha} \leq D_{n,k} \leq \frac{n-1}{(m+\delta)^\alpha}$. Therefore, we have, $\frac{(m+\delta)^\alpha}{n-1} \leq C_{n,k} \leq \frac{(2m+\delta)^\alpha}{n-1}$ w.p.\ 1. \qed

We use Theorem~\ref{thm:4.3} to conclude that $\tilde{D}_{n+1,k}$ and $D_{n+1,k}$ are identically distributed. Therefore, using the previous lemma, we have,
\begin{equation}
\frac{(m+\delta)^\alpha}{n} \leq \tilde{C}^{-1}_{n+1,k} = \tilde{D}_{n+1,k} \leq \frac{(2m+\delta)^\alpha}{n} \label{eq:4.6}
\end{equation}\
as $n \to \infty$.

\begin{prop}\label{prop:4.5}
For all $i \geq 1$, there exists a random sequence $\{c_n\} \sim \Theta(m^2 \log n)$ such that
\[
\frac{\tau_n - \tau_i}{c_n} \to 1 \qquad \text{a.s.}
\]
as $n \to \infty$.
\end{prop}

\emph{Proof.} We observe that the random variables $T_1^{(n+1)}, \ldots, T_m^{(n+1)}$ are independent and that
\begin{align*}
    T_1^{(n+1)} \mid \widetilde{\mathcal{F}}_{n+1,0} &\sim Exp(\tilde{D}_{n+1,0}) \\
    T_2^{(n+1)} \mid \widetilde{\mathcal{F}}_{n+1,1} &\sim Exp(\tilde{D}_{n+1,1}) \\
    &\vdots \\
    T_m^{(n+1)} \mid \widetilde{\mathcal{F}}_{n+1,m-1} &\sim Exp(\tilde{D}_{n+1,m-1})
\end{align*}

We define
\[
b_n = \tilde{C}_{n+1,0} + E[\tilde{C}_{n+1,1} | \widetilde{\mathcal{F}}_{n+1,0}]
    + \cdots + E[\tilde{C}_{n+1,m-1} | \widetilde{\mathcal{F}}_{n+1,0}]
\]

Then, using equation~\eqref{eq:4.6}, we may conclude that,
\begin{equation}
\frac{m(m+\delta)^\alpha}{n} \leq b_n \leq \frac{m(2m+\delta)^\alpha}{n}\label{eq:4.7}
\end{equation}
Now, recall that $\tau_{n+1} - \tau_n = \sum_{k=1}^m T_k^{(n+1)}$. Thus, $\{\tau_{n+1} - \tau_n - b_n, \widetilde{\mathcal{F}}_{n+1,0}\}$ forms a martingale difference sequence. Further, we have, $\sum b_n^2 < \infty$, which implies that $Y_n = \tau_n - \tau_i - (b_i + \cdots + b_n)$ is an $L^2$ bounded martingale. Therefore, by the Martingale Convergence Theorem, we have, $Y_n \to Y$ a.s.\ as $n \to \infty$. We define,
\[
c_n = b_i + \cdots + b_n
\]
Therefore,
\begin{align*}
    \tau_n - \tau_i - c_n &\xrightarrow{\text{a.s.}} Y \\
    \implies \frac{\tau_n - \tau_i}{c_n} &\xrightarrow{\text{a.s.}} 1
\end{align*}
Again, using equation~\eqref{eq:4.7}, it easily follows that $c_n \sim \Theta(m(m+\delta)^\alpha \log n)$. This completes the proof. \qed

\begin{prop}\label{prop:4.6}
For $m = 1$ and $\delta=0$, the sequence $\{c_n\}$ outlined in Proposition~4.7 satisfies $\frac{c_n}{\log n} \to \frac{1}{\lambda^*}$ as $n \to \infty$.
\end{prop}
\emph{Proof.} From Lemma~\ref{lemma:4.1}, we have $n \tilde{C}_n \to \frac{1}{\lambda^*}$ a.s.\ as $n \to \infty$. This observation, along with the form of the sequence $c_n$, help us to conclude the result sought. \qed

\subsection{Proof of Theorem~\ref{thm:ida}\,(1)}
From Theorem~\ref{thm:4.3} we have, for $j=1,\cdots,n$, $d_j(n)$ is distributed identically as $Z_j(\tau_n - \tau_j)$. Also, combining Theorem~\ref{thm:4.4} and Proposition~\ref{prop:4.5}, we have,
\[
\frac{Z_i(\tau_n - \tau_i)}{{c_n}^\frac{1}{1+\alpha}} \to (1+\alpha)^\frac{1}{1+\alpha} \qquad a.s.
\implies
\frac{d_i(n)}{{c_n}^\frac{1}{1+\alpha}} \to (1+\alpha)^\frac{1}{1+\alpha} \qquad a.s.
\]
Finally, we note from Proposition~\ref{prop:4.6}, that $\frac{c_n}{\log n} \to \frac{1}{\lambda^*}$ a.s. This helps us to conclude that the result in consideration.

\subsection{Proof of Theorem~\ref{thm:ida}\,(2)}
We use Theorem~A, \cite{RudTotVal2007} with $\phi(G) = 1(\#\mathrm{children}(\emptyset, G) = k)$ where $\emptyset$ denotes the root of the tree $G$. Then we have,

\[
\lim_{t \to \infty} \frac{|\{x \in \Upsilon(t) : \deg(x, \Upsilon(t)) = k+1\}|}{|\Upsilon(t)|}
= \lambda^* \int_0^\infty \exp(-\lambda^* t) P(\#\mathrm{children}(\emptyset, \Upsilon(t)) = k)\,dt
\]

Now,
\[
    P(\#\mathrm{children}(\emptyset, \Upsilon(t)) = k)
    = P(T_k < t) - P(T_{k+1} < t)
\]
By Fubini's Theorem, we have,
\[
    \lambda^* \int_0^\infty \exp(-\lambda^* t) P(T_k < t) dt = E(e^{-\lambda^* T_k})
\]
$T_k$ is the sum of independent exponentially distributed random variables with parameters $\frac{1}{(1+\delta)^\alpha}, \frac{1}{(2+\delta)^\alpha}, \frac{1}{(3+\delta)^\alpha}, \cdots, \frac{1}{(k+\delta)^\alpha}$; this can be easily calculated. This completes the proof.

\subsection*{4.5. Proof of Corollary~\ref{cor:ida}}
We first define
\[
    \tilde{p}_k = \frac{(k+\delta)^\alpha\lambda^*}{(k+\delta)^\alpha\lambda^* + 1}
    \prod_{i=1}^{k-1} \frac{1}{(i+\delta)^\alpha\lambda^* + 1}, \qquad k = 1, 2, \cdots
\]
(i) We observe that $\lambda^* > 1$ and note that
\[
    \frac{\tilde{p}_{k+1}}{\tilde{p}_k} \leq 1 \quad \forall k \geq 1.
\]
This proves that the mode of the distribution is at $1$.

(ii) We observe that
\begin{align}
    \sum_{k=n}^{\infty} \tilde{p}_k = \prod_{i=1}^{n-1} \frac{1}{1 + (i+\delta)^\alpha\lambda^*}\\
\end{align}

Now,

\begin{align}
    \sum_{n=1}^{\infty}\sum_{k=1}^{\infty}\tilde{p}(k) &= \sum_{k=1}^{\infty}\sum_{n=1}^{k}\tilde{p}(k)\\
    &=\sum_{k=1}^{\infty}k\tilde{p}(k)\\
    &=E[\text{\# of degrees}]\\
    &=2(\because\text{ degrees are preserved})
\end{align}

(iii)
\begin{align}
    \sum_{k=n}^{\infty} \tilde{p}_k &= \prod_{i=1}^{n-1} \frac{1}{1 + (i+\delta)^\alpha\lambda^*}\\
    &= \Theta{\left( \frac{1}{(\lambda^*)^{n-1}(n!)^\alpha} \right)}\\
    &=o\left( \frac{1}{(\lambda^*)^{n-1}(n/2)^{n\alpha/2}} \right)
\end{align}

Clearly $\frac{1}{(\lambda^*)^{n-1}(n/2)^{n\alpha/2}}$ decays at a faster rate than exponential.

\subsection*{4.6. Proof of Theorem~\ref{thm:ida}\,(3)}
We will use the Athreya-Karlin Embedding. Let $\tilde{N}_k(n)$ denote the number of processes with exactly $k$ individuals at time $\tau_n$. \quad \text{Hence, combining Lemma~\ref{lemma:4.1} and Theorem~\ref{thm:ida}\,(2), we have,}
\begin{align}
P(\tilde{d}_v(n+1) - \tilde{d}_v(n) = 1, \tilde{d}_v(n) = k | \widetilde{\mathcal{F}}_n)
&= \tilde{N}_k(n) \cdot \tilde{C}_n \frac{1}{k}\\
&= \frac{\tilde{N}_k(n)}{n} \cdot n \tilde{C}_n \frac{1}{k}\\
&\xrightarrow{a.s.}\prod_{i=1}^k \frac{1}{1 + (i+\delta)^\alpha\lambda^*}
\end{align}

We note that the sequence $\tilde{N}_k(n) \cdot \tilde{C}_n \cdot \frac{1}{k}$ is uniformly bounded and hence uniformly integrable. Taking expectations gives us the result sought.

\subsubsection*{4.6.1. Proof of Theorem~\ref{thm:ida_2}}
From Theorem~\ref{thm:4.3} we have, for $j=1,\cdots,n$, $d_j(n)$ is distributed identically as $Z_j(\tau_n - \tau_j)$. Also, combining Theorem~\ref{thm:4.4} and Proposition~\ref{prop:4.5}, we have,
\begin{align}
\frac{Z_i(\tau_n - \tau_i)}{{c_n}^{\frac{1}{1+\alpha}}} \xrightarrow{a.s.} (1+\alpha)^{\frac{1}{1+\alpha}} \\
\implies
\frac{d_i(n)}{{c_n}^{\frac{1}{1+\alpha}}} \xrightarrow{a.s.} (1+\alpha)^{\frac{1}{1+\alpha}} 
\end{align}
Finally, we note from Proposition~\ref{prop:4.5}, that $c_n \sim \Theta(m(m+\theta)^\alpha \log n)$. This concludes the proof.

%\nocite{*}
\bibliographystyle{plain}
\bibliography{ref}

\end{document}